\documentclass[a4paper, 11pt]{article}
\usepackage{amsmath}
\usepackage{amssymb}
\usepackage{textcomp}
\usepackage{marvosym}
\usepackage{wasysym}
\usepackage{pifont}
\usepackage{graphics}

\title{Handles, Hooks, and Scenarios: A fresh Look at the Collatz Conjecture}
\author{Manfred Tr\"{u}mper}
\date{2006-12-02}

\begin{document}
\maketitle

\noindent \underline{Summary.} An operational approach to the \emph{Collatz Conjecture} is presented. \emph{Scenarios} are defined as strings of characters "s" (for "spike") and "d" (for "down") which symbolize the Collatz operations $(3m+1)/2$ and $m/2$ in a \emph{Collatz Series} connecting two odd integers, the \emph{startnumber} and the \emph{endnumber}. It is shown that a scenario determines uniquely four integers, called \emph{startperiod} $A_M$, \emph{startphase} $B_M$, \emph{endperiod} $A_N$, and \emph{endphase} $B_N$ such that \emph{startnumber} $M = A_M \cdotp k - B_M$, and \emph{endnumber} $N = A_N \cdotp k - B_N$, where k is any natural number. Therefore, any scenario can be realized infinitely many times, with well-defined startnumber and endnumber for each \emph{realization}. It is shown how the periods and phases for a given scenario are calculated. The results are used to prove that any odd (even) number in a Collatz Series is less than 8 (7) ("up"- or "down"-) steps away from an odd integer which is divisible by 3. They are also used for the construction of Collatz Series which exhibit prescribed regular graphics patterns.
Finally, the bearing of the present work on the question of \emph{non-trivial cycles} is discussed.

\section{Introduction.}
The \emph{Collatz Conjecture} has been widely discussed in the mathematical 
literature and numerous authors have contributed to the discussion.
Interesting reviews of the problem have been given by Lagarias (1996) and 
by Wirsching (2000).
\newline 
For this paper a three-fold goal has been set.

First, to present, within a newly developped concept, some new results which hopefully will help to deepen the understanding of the \emph{Collatz Conjecture} and advance future research.

Second, to use key notions which appeal directly to the intuition of the non-mathematician, rather than to the specialist in number theory. Thus we will employ common language terms instead of scientific terms wherever there is no risk of loss of precision. E.g. instead of describing an integer n by the condition $n\equiv 1$ (mod $2$) we will write \textquoteleft n is odd\textquoteright. To a large extent, scientific thinking depends on intuition and not every brain has been trained to think in terms of abstract formalisms.

Third, to emphasize some facts about Collatz Series which have not been accorded sufficient attention in the previous literature. 

The paper is organized as follows. Section 2 introduces the notions of \emph{scenarios}, their \emph{periods}, \emph{phases}, and their \emph{realizations}. Section 3 explains the notion of \emph{hooks} and gives closed form expressions for their periods and phases. Section 4 lays the ground for the practical calculation of periods and phases for a given \emph{scenario}. Sections 5.1 and 5.2 are concerned with an important aspect of \emph{Collatz Backward Series}. Section 5.3 gives an example for the degree of control over Collatz Series which is reached with the results of this paper. With the exception of the outline of a proof for propositions 2.2 and 2.3 and of the proofs for  propositions 3.1, 5.1, and 5.2, no proof for any of the other propositions is given. The proofs were omitted to keep the length of the paper within reasonable bounds. They may be published elsewhere.\\
To facilitate reading the text, the ends of definitions and propositions are marked by the symbol $\Box.$ \\

Let us start by recalling or defining some of the basic terms and notions relevant for the discussion.\\
\textbf{Definition 1.1} \space\space The term \emph{Collatz Series} designates the sequence of  integers which is generated by recursive application of the \emph{Collatz Rules} to a given positive integer $m$:

Determine $3m+1$ if $m$ is odd, and determine $m/2$ if $m$ is even.\newline
The operation $3m+1$ will be denoted by \textquoteleft u\textquoteright\space (for \textquoteleft up\textquoteright), and\\
the operation $m/2$ will be denoted by \textquoteleft d\textquoteright\space (for \textquoteleft down\textquoteright). \space $\Box$ \\

The result of an \textquoteleft u\textquoteright\space operation, $3m+1$, is always an even integer and will, therefore, undergo at least one operation \textquoteleft d\textquoteright.\\

A few more remarks about Collatz-related jargon may be useful.\\
We will not only have to deal with Collatz Series, but also with the series generated by the inverse of the Collatz operations, i.e. the operations $(m-1)/3$ and $2m$. Such a series will be called a \emph{Collatz Backward Series} or, shorter, a \textquoteleft CBS\textquoteright. The ordinary Collatz Series will also be referred to as a \textquoteleft CFS\textquoteright\space (\emph{Collatz Forward Series}).\\

Since the Collatz operations are involving the integers $2$ and $3$, divisibility by $2$ or $3$ will be an important property of the integers encountered. Common language does already describe the divisibility by $2$ in assigning the properties \emph{even} or \emph{odd} to an integer. 
It will be useful as well to establish convenient expressions to describe divisibility by $3$.\\
\textbf{Definition 1.2} \space\space An integer n  will be called \\
\emph{RC0} if $n/3$ has remainder $0$ \space $(\Leftrightarrow n\equiv 0$ \space $(mod 3$)), \\
\emph{RC1} if $n/3$ has remainder $1$ \space $(\Leftrightarrow n\equiv 1$ \space $(mod 3$)), \\
\emph{RC2} if $n/3$ has remainder $2$ \space $(\Leftrightarrow n\equiv 2$ \space $(mod 3$)).
\newline
RC stands for the mathematical term \emph{Residue Class}. \space $\Box$ \newline
E.g. we will say \textquoteleft $6$ is RC0\textquoteright, or \textquoteleft $7$ is RC1\textquoteright, or \textquoteleft $8$ is RC2\textquoteright.  

Here are some simple algebraic properties of these classes which we will need.\\
(a) $2n$ is RC0 iff $n$ is RC0, $2n$ is RC1 iff $n$ is RC2, $2n$ is RC2 iff 
$n$ is RC1.\\
(b) RC0 + RC1 = RC1, RC0 - RC1 = RC2.

With this in mind we can now say that an \textquoteleft up\textquoteright\space operation will always result in an even integer which is RC1. The ensuing \textquoteleft down\textquoteright\space operation will then yield an integer which is RC2. Further division by $2$, if it is possible, will result in an integer which is RC1, and so on. 

For this reason an integer which is RC0 cannot occur in a Collatz Series once an \textquoteleft up\textquoteright\space operation has been made. An odd RC0 integer can occur at most once in a Collatz Series.

Odd RC0 integers are thus singled out among the other integers in a Collatz Series.

Also, since an odd RC0 integer can only be reached (from even RC0 integers) by  \textquoteleft down\textquoteright\space operations, the backward extension of a Collatz Series from there is unique. This means that an odd RC0 integer determines a Collatz Series uniquely, both in forward and in backward direction.

We have to recall here that starting with an RC1 or RC2 integer, the Collatz Backward Series (CBS) is \emph{not} unique, i.e. for each such integer there is a backward tree with an infinte number of branches.
However, there is only one way to go backwards from a RC0 integer, namely by repeated doubling.

In view of the above it is appropriate to give the odd RC0 integers a proper name.\\
\textbf{Definition 1.3} \space\space  An integer which is odd and RC0 (i.e. divisible by $3$ with remainder $0$) will be called a \emph {handle}. \space $\Box$
\newline
The term handle was chosen since both the Collatz forward series (CFS) and the Collatz backward series (CBS) are determined uniquely by such an integer. They are essential for the notion of the \emph{Complete Collatz Series}. It comes down from infinity through divisions by $2$, with RC0 integers of the form $2^{n}H$ ($H = 3M$, $M$ an arbitrary odd integer), until it reaches the handle H. From thereon the series is made up of integers RC1 or RC2.

We shall show later (Section 5.1) that every odd integer can be reached in less than $8$ steps by a Collatz Series starting from a handle. 
\section{Scenarios and their Realization.}

Since an operation \textquoteleft u\textquoteright\space is always followed by an operation \textquoteleft d\textquoteright, the combined operation \textquoteleft ud\textquoteright\space will be denoted by \textquoteleft s\textquoteright\space (for \emph{spike}, the graphical appearance of the operation).

It is clear that any Collatz Series can be described in terms of \textquoteleft spikes\textquoteright\space and \textquoteleft downs\textquoteright.\newline
\textbf{Definition 2.1} \space A sequence of Collatz operations connecting 2 odd integers will be called a \emph{scenario}. \space $\Box$  
\newline
A scenario is thus described by a word, beginning with \textquoteleft s\textquoteright\space  and followed by any combination of the characters \textquoteleft s\textquoteright\space or \textquoteleft d\textquoteright. \newline
\textbf{Definition 2.2} \space A series of integers which is described by a specified scenario will be called a \emph{realization of the scenario}. \space $\Box$ \newline
As any scenario starts with an \textquoteleft s\textquoteright, the realization of a scenario starts with an odd integer.

A scenario literally spells out what happens between a startnumber $M$ 
and an endnumber $N$.

For a scenario with $\sigma$ operations \textquoteleft s\textquoteright\space and $\delta$ operations \textquoteleft d\textquoteright\space the general form is 
\begin{center} 
\fontsize{14pt}{18pt} $sd^{\delta_{1}}sd^{\delta_{2}}...sd^{\delta_{\alpha}}...sd^{\delta_{\sigma}}$ 
\normalsize
\end{center}
where $\delta_{1}+\delta_{2}$ +...+ $\delta_{\sigma} = \delta $. 
The exponents $\delta_{\alpha}$ ($\alpha = 1, 2, ...\sigma$) are describing the multiplicity of the respective operation \textquoteleft d\textquoteright.\\

We are now arriving at a crucial point in our discussion. The simplest, almost trivial question about Collatz Series is to ask, \emph{'Given a startnumber, what is the series?'} This question is, of course, the same as asking \emph{'Given a startnumber, what is the scenario attached to it?'} Many authors of articles about the Collatz Conjecture have posed and answered this question.\\

Here we will turn it around and ask the simple, but non-trivial question \begin{center}\emph{\textquoteleft Given a scenario, what are the startnumbers?\textquoteright}\end{center}

Our task will be to find the startnumbers which are compatible with a given scenario. In fact, we want to determine all realizations of a scenario.\newline\\\\
\textbf{Proposition 2.1 (Periodicity of a Realization)} \space\space 
A scenario $S$ uniquely determines \newline
- a pair of integers $A_M$ and $B_M$, with $A_M$ even, $B_M$ odd,\newline
- a pair of integers $A_N$ and $B_N$, with $A_N$ even, $B_N$ odd, \newline
such that for any $k = 1, 2, 3,$ ...... \newline
$M_k = A_M\cdot k - B_M$  is the startnumber, and\newline
$N_k = A_N\cdot k\space - B_N$ is the endnumber \newline
of the Collatz sequence described by the scenario. \space $\Box$ \\\\
\textbf{Corollary 2.1} \space The RC-property of the startnumber $M_k$ and the endnumber $N_k$ is not changed under the replacement $k$\textrightarrow $k+3p$ where p is any (positive) integer.$\Box$\\

\noindent\textbf{Definition 2.3} \space\space The \emph{Collatz Series} starting with $M_k$ and ending with $N_k$ will be referred to as the \emph{kth realization} of the scenario.
The integer $A_M$  ($A_N$) will be called the \emph{startperiod} (\emph{endperiod}).
The integer $B_M$  ($B_N$) will be called the \emph{startphase} (\emph{endphase}). \space  $\Box$ \newline\\
\textbf{Proposition 2.2 (Startperiod)} \space\space Let $\sigma$ and $\delta$ denote, respectively, the numbers of \emph{spikes (\textquoteleft s\textquoteright)} and \emph{downs (\textquoteleft d\textquoteright)} in a scenario.
The startperiod depends on the total number of operations \textquoteleft s\textquoteright\space and \textquoteleft d\textquoteright\space and is given by $A_M = 2\cdot 2^{\sigma+\delta}$. \space $\Box$ \newline\\
\textbf{Proposition 2.3 (Endperiod)} \space\space The endperiod depends only on the number of spikes in the scenario 
and is given by $A_N = 2\cdot 3^{\sigma}$. \space $\Box$ \newline 

\underline{Outline of a proof.} 
Consider the simplest and most basic scenario which is just \textquoteleft s\textquoteright.
Apply the operation $(3M + 1)/2$ to an unspecified odd startnumber $M = 2k - 1$, where $k$ may be any positive integer. The result is $3k - 1$ which is even or odd depending on $k$. To get an odd endnumber, re-scale $k$ to $2k$. This yields $N_{k} = 6k - 1$. To get the correct values of the startnumber, $M$ has to be re-scaled accordingly and then becomes $M_{k} = 4k - 1$. To verify that these expressions are realizations of the scenario \textquoteleft s\textquoteright\space for all integers $k$, just apply the Collatz rules to $M_{k}$.\space$\blacksquare$

By similar reasoning one can derive the realizations for the next simplest scenarios, with $k = 1, 2, 3...$, \\
\textquoteleft sd\textquoteright\space : Startnumbers are $M_{k} = 8k - 7$, endnumbers are $N_{k} = 6k - 5$,   \newline
\textquoteleft ss\textquoteright\space : Startnumbers are $M_{k} = 8k - 1$, endnumbers are $N_{k} = 18k - 1$.
\\
Again, these expressions may be verified by applying the Collatz rules to the startnumbers.\newline

Unlike the periods, the phases are not given by simple expressions. The reason is that they depend on other properties of a scenario, such as the order in which the \textquoteleft s\textquoteright\space and \textquoteleft d\textquoteright\space appear. However, the following can be said.\newline\newline
\textbf{Proposition 2.4}\space\space Both startphase and endphase are odd integers.\newline 
The startphases satisfy the condition $0$ \textless $B_M$ \textless $A_M$. \newline 
The endphases satisfy the conditions  $0$ \textless $B_N$ \textless $A_N$ and, in addition, they are RC1 or RC2 (and not RC0). \space $\Box$ \newline \newline Further remarks about endphases will be found in section 4.
These conditions are sufficient to assure that, in all realizations of the scenario, startnumbers and endnumbers are odd integers and, in addition, that the endnumbers are not divisible by $3$.\newline 
\section{Hooks.}
We will now take a look at the building blocks of any Collatz Series.\\
\textbf{Definition 3.1}\space\space A $\delta$-hook is a scenario of the type $sd^\delta$. The exponential notation indicates that the \textquoteleft down\textquoteright\space operation is performed $\delta$ times.\space$\Box$ \newline

Hooks got their name from the way they appear in graphics. Since all scenarios are sequences of hooks hooked up to other hooks, this notion is basic for the discussion of the Collatz conjecture. 

Period and phase for startnumber and endnumber of a hook are relatively easy to calculate. 
It turns out that the periods don't, but the phases do depend on the parity of $\delta$.\newline\newline
\textbf{Proposition 3.1} \space\space For a $\delta$-hook the following holds.\newline
\emph {Periods:}\space\space Startperiod:\space $A_M$ = $2^{\delta+2}$,\space\space Endperiod:	$A_N = 6$,\newline
in accordance with propositions 2.2 and 2.3.\newline
\emph {Phases:}\newline
\underline{$\delta$ even:}\space	Startphase:\space\space$B_M = (2^{\delta+1} + 1)/3$ ,\space\space \space \space \space \space Endphase:\space\space$B_N = 1$. \\
\underline{$\delta$ odd:}\space\space	Startphase:\space\space$B_M = (5\cdot 2^{\delta+1} +1)/3$ ,\space\space Endphase:\space\space$B_N = 5$. $\Box$ \newline\\
\underline {Proof.}\space The startnumbers $M_k$ and endnumbers $N_k$ of a $\delta$-hook are given by\newline
\underline{$\delta$ even:}		$M_k = 2^{\delta+2}\cdot k - (2^{\delta+1} + 1)/3$ ,\space \space \space \space \space	$N_k = 6 k - 1$ \newline
\underline{$\delta$ odd:} \space	$M_k = 2^{\delta+2}\cdot k - (5\cdot 2^{\delta+1} + 1)/3$ , 	$N_k = 6 k - 5$ \newline
with $k = 1, 2, 3,...$ \newline
Using these expressions we obtain by straightforward calculation\\
$\delta$ even: $(3M_k+1)/2 = 2^{\delta}(6k-1)$,\\
$\delta$ odd:\space\space $(3M_k+1)/2=2^{\delta}(6k-5)$.\space$\blacksquare$\\

It should be noted that the endphase of a hook with $\delta$ even (odd) is RC1 (RC2). Therefore, since the endperiod of a hook is always $6$, the endnumber of a hook with $\delta$ even (odd) is RC2 (RC1).\newline
In view of the relevance of the parity of $\delta$ we shall call a hook with $\delta$ even (odd) an \emph{even (odd) hook}.
Since the endnumber of any scenario and the last hook in it have the same RC-property, we have  \newline
\textbf{Proposition 3.2} \space\space The endnumber of a scenario is RC2 (RC1) if the last hook is even (odd).\space $\Box$ \newline

This simple fact is also evident from basics since the last (in fact, any) \textquoteleft s\textquoteright\space operation in a scenario sets the RC-property of the integer reached to $2$.\\\\
Values for the periods and phases of the first 16 hooks are given in the table below.\\\\
\begin{quote}
\begin{tabular}{|r|r|r|c|c|}
\hline
$\delta$&Startperiod&Startphase&Endperiod&Endphase\\
\hline
0&4&1&6&1\\
1&8&7&6&5\\
2&16&3&6&1\\
3&32&27&6&5\\
4&64&11&6&1\\
5&128&107&6&5\\
6&256&43&6&1\\
7&512&427&6&5\\
8&1.024&171&6&1\\
9&2.048&1.707&6&5\\
10&4.096&683&6&1\\
11&8.192&6.827&6&5\\
12&16.384&2.731&6&1\\
13&32.768&27.307&6&5\\
14&65.536&10.923&6&1\\
15&131.072&109.227&6&5\\
\hline
\end{tabular}
\end{quote}
\begin{center}
\textbf {The first 16 hooks with their periods and phases.}\\
\end{center}

Figure 1 shows the graphics for various combinations of $\delta$-hooks with $\delta = 0, 1, 2$. Each hook is pictured as the straight line connecting the startnumber with the endnumber. Since the intermediate even integers are not shown, the \textquoteleft hook appearence\textquoteright\space gets lost. Each column shows a number of realizations for the type of hook indicated below the horizontal axis. You have to imagine that each of these columns reaches upwards towards infinity! The number of realizations shown in the Figure is only limited by the space available on paper.\\ The graphics shows how different hooks will connect to form new scenarios.\\1. The scenario \textquoteleft ssd\textquoteright\space which is made up of hooks from the first two columns, is shown with its first 10 realizations. This scenario has startnumbers $M_k=16k-5$ and endnumbers $N_k=18k-5$.\\
2. The scenario \textquoteleft sdsdd\textquoteright\space which is made up of hooks from the second and third column, is shown with its first 4 realizations. This scenario has startnumbers $M_k=64k-47$ and endnumbers $N_k=18k-13$.\\3. The scenario \textquoteleft sdds\textquoteright\space which is made up of hooks from the third and fourth column, is shown with its first 5 realizations. This scenario has startnumbers $M_k=32k-3$ and endnumbers $N_k=18k-1$.\\4. The scenario \textquoteleft ssdd\textquoteright\space which is made up of hooks from the fourth and fifth column, is shown with its first 2 realizations. This scenario has startnumbers $M_k=32k-13$ and endnumbers $N_k=18k-7$.\\5. The scenario \textquoteleft ssdsdd\textquoteright\space which is made up of hooks from the first three columns, is shown with its first 2 realizations. This scenario has startnumbers $M_k=128k-117$ and endnumbers $N_k=54k-49$.\\6. The scenario \textquoteleft sdsdds\textquoteright\space which is made up of hooks from columns 2 through 4, is shown with its first 2 realizations. This scenario has startnumbers $M_k=128k-47$ and endnumbers $N_k=54k-19$.\\7. The scenario \textquoteleft ssdsdds\textquoteright\space which is made up of hooks from columns 1 through 4, is shown with its first 2 realizations. This scenario has startnumbers $M_k=256k-117$ and endnumbers $N_k=162k-73$.\\\\The trivial cycle is represented by the hook \textquoteleft sd\textquoteright\space in first realization. It is shown by the horizontal line at the bottom of the second column. Now, if there would be a non-trivial cycle it would be shown in a similar diagram, with many more hooks of various kind arranged in various order, and it would have a startnumber to the left at the same hight as the endnumber to the right.\\
\section{Calculation of Periods and Phases by Iteration.}
Here we will describe how the periods and phases are changed when one operation is appended to a scenario.
 
Appending an \textquoteleft s\textquoteright\space is easy because this operation is acting directly on the odd endnumber of the given scenario. 

Adding a \textquoteleft d\textquoteright\space is a bit more tricky since this operation does not act on the endnumber of the given scenario. \newline
But in any case we start from a scenario $\mathcal{S}$ with realizations 

	$M_k = A_M\cdot k - B_M$ and $N_k = A_N\cdot k - B_N$  , 
where $k = 1, 2, 3, ...$ \newline
The enlarged scenario $\mathcal{S'}$ has realizations with \newline 
startnumbers $M'_k = A'_M\cdot k - B'_M$ and endnumbers $N'_k = A'_N\cdot k - B'_N$  , where $k = 1, 2, 3, ...$ \newline\newline
\textbf{Proposition 4.1 (Appending an \textquoteleft s\textquoteright)} \space Let $\mathcal{S'}$ = $\mathcal{S}$s be the scenario with one operation \textquoteleft s\textquoteright\space added at the end of $\mathcal{S}$.\newline
The periods and phases of $\mathcal{S'}$ are expressed by those of $\mathcal{S}$ as follows.\newline\newline
\underline{\textonehalf$(3B_N - 1)$ odd:}\smallskip\\
	$A'_M = 2A_M$,\space  $B'_M = B_M$ , \newline
	$A'_N = 3A_N$,\space\space  $B'_N$ = \textonehalf$(3B_N - 1)$. \newline\newline
\underline{\textonehalf$(3B_N - 1)$ even:}\smallskip\newline
	$A'_M = 2A_M$,\space  $B'_M = B_M + A_M$ ,\newline
	$A'_N = 3A_N$,\space\space  $B'_N$ = \textonehalf$(3B_N + 3A_N - 1)$. $\Box$  \newline\newline It is worth noting that in the case where \textonehalf$(3B_N - 1)$ is \emph{odd}, the startphase remains unchanged. Simple reasoning shows that this is true for any $B_N$ which can be written as $B_N = 4j-3$, with any $j$ that is RC1 or RC2 (not RC0!).\\Likewise, in the case where \textonehalf$(3B_N - 1)$ is \emph{even}, the new startphase is obtained by simply adding the startperiod to the startphase of the given scenario. The condition is satisfied by any $B_N$ which can be written as $B_N = 4j-1$, with any $j$ that is RC0 or RC2 (not RC1!).\\The restrictions on j are necessary to ensure that $B_N$, and by implication $N_k = A_{N}k-B_N$, is not a handle.\\\\
\textbf{Proposition 4.2 (Appending a \textquoteleft d\textquoteright)} \space Let $\mathcal{S'}$ = $\mathcal{S}$d be the scenario with one operation \textquoteleft d\textquoteright\space added at the end of $\mathcal{S}$.\newline
The periods and phases of $\mathcal{S'}$ are expressed by those of $\mathcal{S}$ as follows.\newline\newline
\underline{\textonehalf$(B_N +$ \textonehalf $A_N)$ odd:}\smallskip\newline
	$A'_M = 2A_M$,\space\space\space  $B'_M = B_M + $\textonehalf $A_M$ ,\newline
	$A'_N =  A_N$,\space\space\space\space $B'_N$ = \textonehalf$(B_N + $\textonehalf $A_N)$ .\newline\newline
\underline{\textonehalf$(B_N + $\textonehalf $A_N)$ even:} \smallskip \\ 
	$A'_M = 2 A_M,\space\space  B'_M = [B_M + 3A_M/2 ]$ mod $A_M$,\newline
	$A'_N =  A_N$,\space\space\space\space  $B'_N$ =  [\textonehalf$(B_N + 3A_N/2)]$ mod $A_N$. $\Box$ \newline 

The above formulas may be used (and have been used by the author) to write a program which calculates all periods and phases as well as the startnumber and endnumber for the $k$th realization of any given scenario, see next section.\newline
\section{Applications and Examples.}
Given that handles are playing such a unique role in Collatz Series it is interesting to ask how far they are away from ordinary integers, i.e. integers which are even or odd, RC1 or RC2. By \textquoteleft how far away\textquoteright\space we mean the number of steps (up or down) by which the integer can be reached via a Collatz Series.\\
For reasons which will become clear in the following it turns out to be advantageous to consider separately odd and even target integers.\\\\
\textbf{5.1 Hooking up \emph{odd} Integers to Handles.}\newline
\textbf{Proposition 5.1}\space\space\ Any odd integer is located on a Collatz Series within less than $8$ steps from a handle. (One step is either an \textquoteleft up\textquoteright\space or a \textquoteleft down\textquoteright\space operation.) $\Box$ \newline\\
\underline{Proof.} If the integer in question is odd and RC0, nothing has to be shown because it is itself a handle. We will show that hooks with $\delta$ \textless \space $6$ can provide the link between handles and any odd (none RC0) integer. 
These hooks are listed in the first column of the table below, followed by the number of (up or down) steps in the second column. In the fifth column we have entered the smallest values of k for which the corresponding startnumber $M_k$ is a handle. These startnumbers are listed in the sixth column (headed \textquoteleft $M$\textquoteright). The last column shows the corresponding endnumbers.
\begin{quote}
\begin{tabular}{|c|c|r|c|r|r|r|}
\hline
Hook& Steps& Startnumber $M_k$& Endnumber $N_k$& $k$& $M$& $N$\\
\hline
$s$&2&$4k-1$&$6k-1$&1&3&5\\
$sd$& 3& $8k - 7$& $6k - 5$& 2& 9& 7\\
$sd^2$& 4& $16k - 3$& $6k - 1$& 3& 45& 17\\
$sd^3$& 5& $32k - 27$& $6k - 5$& 3& 69& 13\\
$sd^4$& 6& $64k - 11$& $6k - 1$& 2& 117& 11\\
$sd^5$& 7& $128k - 107$& $6k - 5$& 1& 21& 1\\
\hline
\end{tabular}
\end{quote}

The endnumbers (last column) comprise all odd non-RC0 integers which are in the interval (0, 18).
It can be verified by direct calculation that for each scenario the target is reached.

Now replace k (in col. 5, table above) by k + 3p where p is any positive integer. 
The startnumbers will become larger, but they will remain handles.
The endnumbers will just be increased by 18p.\\
For $p=0$ the k-values, startnumbers, and endnumbers are shown in the table above (col. 5, 6, 7).\\
For $p=1$ the resulting k-values, startnumbers, and endnumbers are shown in the table below (columns 5, 6, 7) and for $p=2$ they are shown in columns 8, 9, 10.\newline

\begin{quote}\space 
\begin{tabular}{|c|r|r|r|c|r|r|c|r|r|}
\hline
Hook&$k+3p$&$M_p$&$N_p$ &$k$&$M$&$N$&$k$&$M$&$N$\\
\hline
$s$&$1+3p$&$3+12p$&$5+18p$&4&15&23&7&27&41\\
$sd$&$2+3p$&$9+24p$&$7+18p$&5&33&25&8&57&43\\
$sd^2$&$3+3p$&$45+48p$&$17+18p$&6&93&35&9&141&53\\
$sd^3$&$3+3p$&$69+96p$&$13+18p$&6&165&31&9&261&49\\
$sd^4$&$2+3p$&$117+192p$&$11+18p$&5&309&29&8&501&47\\
$sd^5$&$1+3p$&$21+384p$&$1+18p$&4&405&19&7&789&37\\
\hline
\end{tabular}
\end{quote}

If we let p run through all natural numbers 0, 1, 2, 3, ....
all odd non-RC0 integers will occur as endnumbers in one of the hooks which begins with a handle.\space $\blacksquare$ \newline\newline
The result does not mean, that a hook always provides the shortest way to  connect an odd integer to a handle. For instance, the integer $19$ is reached, via a hook, from the handle $405$ in $7$ steps as shown in the table above. But it is reached in $6$ steps from the handle $33$ (via the scenario sdsd). \newline\newline
\textbf{5.2 Hooking up \emph{even} Integers to Handles.}\newline
Proposition 5.1 applies to all \emph{odd} integers which are not handles, i.e. the integers beginning 1,5,7,11,13,17,19,23,25,29,31,35,37,41,43,47,49,53, etc.
Naturally, the question arises if a similar result holds for \emph{even} integers.\\ 
\textbf{Proposition 5.2}\space Any even integer (not RC0) is located on a Collatz Series within less than $7$ steps from a handle. (One step is either an \textquoteleft up\textquoteright\space or a \textquoteleft down\textquoteright\space operation.) $\Box$ \newline
\underline{Proof.}  
First of all it is clear that the doubles of the odd integers above will be found within less than $7$ steps from a handle. These are the integers 2,10,14,22,26,34,38,46,50,58,62,70,74,82,86,94,98,106, etc.\\
But how about the even integers in between, beginnning with \\ 4,8,16,20,28,32,40,44,52,56,64,68,76,80,88,92,100,104, etc.?\\
Note that the above series is made up of the $2$ sub-series\\
$12j-8$ and $12j-4$ where $j=1,2,3...$.\\ It is easily verified that the integers in the first sub-series are hooked up to handles in the following way (with p taking the values 1,2,3...).\\
For $j=3p$ the endnumber $36p-8$ is obtained from the handle $12p-3$ in one step (\textquoteleft up\textquoteright).\\
For $j=3p-1$ the endnumber $36p-20$ is obtained from the handle $48p-27$ in three steps (by the 1-hook \textquoteleft sd\textquoteright).\\
For $j=3p-2$ the endnumber $36p-32$ is obtained from the handle $192p-171$ in five steps (by the 3-hook \textquoteleft sddd\textquoteright).\\
Similarly, the integers in the second sub-series are hooked up to handles
as follows.\\
For j = $3p$ the endnumber $36p-4$ is obtained from the handle $24p-3$ in two steps (by the 0-hook \textquoteleft s\textquoteright).\\
For $j=3p-2$ the endnumber $36p-28$ is obtained from the handle $96p-75$ in four steps (by the 2-hook \textquoteleft sdd\textquoteright).\\
For $j=3p-1$ the endnumber $36p-16$ is obtained from the handle $384j-171$ in six steps (by the 4-hook \textquoteleft sdddd\textquoteright).$\blacksquare$ \\

What is the significance of the result? \newline
It shows that handles are lurking around the corner near all integers of a Collatz Series. And, in the unlikely case that a non-trivial cycle should exist, it would be surrounded by a \emph{halo} of close-by handles. People who have dealt with Collatz Backward Series (CBS) may have noticed that non-trivial CBS (i.e. those which are not obtained by merely \emph{doubling} RC1 or RC2 integers) will sooner or later hit a handle. And, recalling that one in 3 integers is a handle, one realizes how difficult it would be  for the CBS on a cycle to miss all the handles that are in its way.  \newline

But, what really counts here is the way of treating the problem: The periodicity of scenarios is used to transport some property about Collatz Series from one finite interval (the period of a scenario) to all intervalls, i.e. to infinity. The method is like a \textquoteleft mathematical telescope\textquoteright\space allowing to make statements about numbers which are so large that there would not be enough available paper on earth to write them down, not even in fine print!\newline\newline
\textbf{5.3 Collatz Series Sections without \emph{Hailstones}.}\newline
Some of the publications on the Collatz conjecture have invoked the notion of \emph{hailstones} to describe the seemingly random order in which integers are appearing in a Collatz Series. \newline
To exemplify the benefits which can be drawn from propositions 4.1 and 4.2 we want to produce a startnumber which leads to the section of a Collatz Series that has a perfectly regular pattern. For example, let us look at the scenario $\mathcal{S}=(s^7d^4)^9$. The exponents indicate the multiplicity of the respective operation. So, in this scenario we have in total $7 \cdot 9 = 63$ operations \textquoteleft s\textquoteright\space and $4 \cdot 9 = 36$ operations \textquoteleft d\textquoteright. In terms of the traditional up- and down operations these figures amont to 162 steps.\newline
Beginning with the scenario \textquoteleft s\textquoteright, with startnumbers $M_{k} = 4k - 1$ and endnumbers  $N_{k} = 6k - 1$, we apply the results of the previous section 4 and perform the $98$ remaining iterations to get the periods and phases of our scenario. The result is the following.\\\\
$Startperiod = 1267650600228229401496703205376$\newline
$Startphase = 1039655887956965120651972413057$\newline
$Endperiod = 2289122546861674989771899392854$\newline
$Endphase = 1877409858577201070748176480485$.\newline 
Each of these 4 integers has the order $10^{30}$.\newline\newline
The first 3 realizations of the scenario are as follows.\newline
\underline{$1^{st}$ Realization:}\\
$Startnumber = 227994712271264280844730792319$,\newline
$Endnumber = 411712688284473919023722912369$.\newline
Each of these 2 integers is of the order $10^{29}$.\newline
\underline{$2^{nd}$ Realization:}\\
$Startnumber = 1495645312499493682341433997695$,\newline
$Endnumber = 2700835235146148908795622305223$.\newline
\underline{$3^{rd}$ Realization:}\\
$Startnumber = 2763295912727723083838137203071$,\newline
$Endnumber = 4989957782007823898567521698077$.\newline
Each of these 4 integers is of the order $10^{30}$.\newline

Using a suitable software for handling large integers and the iteration formulas of section 4, all of the above integers were calculated in just the fraction of a second.\newline\\
\emph{All-Number graphics and Odd-Number graphics.} When it comes to the graphical representation of Collatz Series, we have to make the choice between two options. We may plot either \emph{all} integers which are generated by up- or down operations, creating a \emph{AN graphics} (\emph{AN} stands for \textquoteleft All Numbers\textquoteright), or we may plot only the \emph{odd} integers which occur in a CFS, creating a \emph{ON graphics} (\emph{ON} stands for \textquoteleft Odd Numbers\textquoteright). The latter has the clear advantage of displaying much less jitter than we would see in the AN graphics.\\\\
Figure 2 pictures, as ON-graphics, the first 3 realizations of the scenario. The curves show a nearly horizontal and slightly increasing regular zigzag pattern with the endnumbers less than twice the corresponding startnumbers. Once the scenario is terminated, the habitual, chaotic looking shape of the series is back. Figure 3 shows, also as ON-graphics, the section of the Collatz Series presentating the first 3 realizations of the scenario. In Figure 4, the scenario sections of the Collatz Series are shown in full detail, i.e. with \emph{all} numbers pictured (AN-graphics).  \newline Please note also the difference in counting steps in ON-graphics and AN-graphics. In ON-graphics (Figures 2, 3) we just count the odd integers (whose number equals the number of \textquoteleft s\textquoteright\space operations in any section of a Collatz Series). In AN-graphics (Figure 4) we count all integers, i.e. for a scenario with $\sigma$ \textquoteleft s\textquoteright-operations and $\delta$ \textquoteleft d\textquoteright-operations, the number of steps counted is $2\sigma + \delta$.\\

The result shows that the integers appearing in a Collatz Series are not always falling like \textquoteleft hailstones\textquoteright. Rather they are what we want them to be, within the limits of what is allowed by the Collatz rules. We could construct Collatz Series with much longer scenarios than the one considered here, expressing any regular pattern. Most Collatz Series are looking like random patterns because the startnumbers were chosen at random!\newline
In the above example, the scenario was derived from the requirement that, 
in first realization, the endnumber $N$ of the subscenario $s^7d^4$ should be close to the startnumber $M$, i.e. satisfy the condition $\rho = $\textbar$(M-N)/N$\textbar\space $\ll 1$. \\For our scenario, $(s^7d^4)^9$, we have $\rho \simeq 0.45$, but for the subscenario $s^7d^4$ it holds $\rho  \simeq 0.064$.\\\newpage
\section{Epilogue and Outlook.}
\textbf{Epilogue.} After the present work on scenarios had been done, G. J. Wirsching pointed out that it might be related to that of Riho Terras. The paper of Terras appeared 1976 in Acta Arithmetica and can be found via a link at the web site http://matwbn.icm.edu.pl/. 
Terras work (which is not quite easy to understand by somebody who is not specialized in mathematics) is indeed related to the present work. In the following, that relation will be described, employing the terminology developped in this paper. 
The \emph{encoding vector} defined by Terras is closely related, but logically not the same as the \emph{scenario} we defined here. Both notions have in common that they characterize a sequence of Collatz operations. If you consider a \emph{scenario} and replace in it each operation \textquoteleft s\textquoteright\space by \textquoteleft 1\textquoteright\space and each operation \textquoteleft d\textquoteright\space by \textquoteleft 0\textquoteright, you get the corresponding \emph{encoding vector}. (We remark in passing that the \emph{encoding vector} is just a sequence of Boolean numbers, but does not really qualify as a \emph{vector} since this term has a specific meaning in geometry). The difference of the two notions lies in the situations they cover. The \emph{scenario} is defined such that the corresponding Collatz sequence will start and end with an \emph{odd} integer. In contrast, the \emph{encoding vector} corresponds to a Collatz sequence which may end in an even or an odd integer. For this reason, Terras finds for the period of his \emph{encoding representation} (which is his startnumber) half the value of that which we have determined for the \emph{scenario}. In other words, there are twice as many (numerical) realizations for an \emph{encoding vector} than for a \emph{scenario}. 
Furthermore, since the endnumbers may be even, the Terras \emph{encoding vectors} cannot always be connected with each other. In contrast, \emph{scenarios} can always be connected since they always end with an odd integer.
Nevertheless, the result of Terras (his Theorem 1.2) which he termed himself \emph{a remarkable periodicity phenomenon} stands out as a singular achievement.\\\\
\textbf{Outlook.} Where should future work go from here?\\ The first important goal should be to disprove the existence of non-trivial cycles. The trivial cycle (the first realization of the 1-hook \textquoteleft sd\textquoteright, see Figure 1) has only one odd integer in it. If a non-trivial cycle exists it would have a great number (at least about $10^7$) of odd integers in it. Every odd integer in a cycle could be taken as the start- and endnumber of a scenario. If the number of odd integers is $\sigma$, there would be $\sigma$ scenarios to take into consideration. All these scenarios would have the same startperiod and the same endperiod. Their startphases would be different from each other and so would be their endphases. There would be $\sigma$ equations to be satisfied, stating for each scenario the condition that startnumber and endnumber are the same. If these equations would be incompatible, the question of non-trivial cycles would be settled.\\\\
Another question of interest for the Collatz problem would be the composition (i.e. the concatenation) of scenarios. How to determine the phases of the scenario $\mathcal{S} = \mathcal{S}_1 \mathcal{S}_2$ in terms of properties of the scenarios $\mathcal{S}_1$ and $\mathcal{S}_2$? As to the \emph{periods} there is no problem, since from propositions 2.2 and 2.3 we simply find the startperiod $A_M =$ \textonehalf $A^1_M\cdot A^2_M$ and the endperiod $A_N =$ \textonehalf $A^1_N\cdot A^2_N$. However, a much harder nut to crack would be to express the \emph{phases} of $\mathcal{S}$ by the periods and phases of $\mathcal{S}_1$ and $\mathcal{S}_2$.\\
As a first step one should try to have a look at hooks. Since hooks are the building blocks of scenarios it would be useful to have rules for their composition. First calculations on the composition of a j-hook and a k-hook indicate that the result will depend, among other factors, on the RC-property of k. \\\\
Still another interesting question would be about \emph{flat CBS}. These are obtained by performing, each time an odd RC1 or RC2 integer is reached, only the minimal number of doublings necessary to get an even RC1 integer. That minimal number of doublings would be equal to one if the odd integer is RC2 and it would be equal two if the odd integer is RC1. Consider one fixed integer and the set of CBS going back from there. Among all CBS emanating from that integer, the flat CBS and the hooks are the two opposite  extremes. While hooks  represent the steepest increase in a CBS, flat CBS are the other extreme. Can an upper limit be posed on the number of steps necessary to reach a handle via a \emph{flat CBS}?\\\\
\textbf{References.}\\
J.C. Lagarias, \textquotedblleft The 3x+1 Problem and its Generalizations\textquotedblright, Amer.Math.Monthly \textbf{92}, pp. 3-23 (1985)\\ 
http://oldweb.cecm.sfu.ca/organics/papers/lagarias/paper/html/paper.html, 1996\\
R. Terras, \textquotedblleft A stopping time problem on the positive integers\textquotedblright, Acta Arithm. \textbf{30}, (3), 241-252, 1976 \\
G.J. Wirsching, \textquotedblleft Das Collatz-Problem\textquotedblright, Lexikon der Mathematik,\\ Spektrum Akademischer Verlag, Band 1, 335-340, 2000 \\\\\\
Author's Address: 1, chemin du Peiroulet, F-30700 Uz\a`{e}s\\
e-mail: matrumper@tele2.fr
\end{document}